\documentclass[12pt]{article}
 \usepackage{epsfig}
\usepackage{graphics}

\usepackage{color}
\usepackage{epsfig,cmmib57}
\usepackage{amssymb,amsmath,exscale}
 \numberwithin{equation}{section}

\newcommand{\ZSUno}{\sum _{n=1}^{+\ZIN}}
\newcommand{\ZOMq}{\Omega}

\newcommand{\zben}{\beta_n}

\newcommand{\intp}{\int_0^{\pi}}

\newcommand{\intt}{\int_0^t}
\newcommand{\ints}{\int_0^s}

\newtheorem{Theorem}{Theorem}

\newtheorem{Lemma}[Theorem]{Lemma}

\newtheorem{Remark}[Theorem]{Remark}
\newcommand{\zdiaform}{\mbox{~~\zdia}}

\newcommand{\zaa}{\alpha}

\newcommand{\ZDE}{\delta}
\newcommand{\zt}{\tau}
\newcommand{\zdia}{~~\rule{1mm}{2mm}\par\medskip}
\newcommand{\zthe}{\theta}

\newcommand{\ZLA}{\label}
\newcommand{\ZIN}{\infty}
\newcommand{\zProof}{{\noindent\bf\underbar{Proof}.}\ }

\newcommand{\ZBI}{\bibitem}
\newcommand{\ZD}{\;\mbox{\rm d}}

\renewcommand{\zthe}{w}

\author{
L. Pandolfi\thanks{Dipartimento di Scienze Matematiche ``Giuseppe Luigi Lagrange'', Politecnico di Torino, Corso Duca degli Abruzzi 24, 10129 Torino, Italy (luciano.pandolfi@polito.it)}
}
\title{Identification of a relaxation kernel using  two boundary measures\thanks{
This papers fits into the research program of the GNAMPA-INDAM and has been written in the framework of the   ``Groupement de Recherche en Contr\^ole des EDP entre la France et l'Italie (CONEDP-CNRS)''.}}

\begin{document}

\maketitle
{\bf\underline{Abstract}:}  
We consider a distributed system with persistent memory
of a type which is often encountered in viscoelasticity or in the study of diffusion processes with memory. The \emph{relaxation kernel,} i.e. the kernel of the memory term, is scarcely known from first principles, and it has to be inferred from experiments taken on samples of the  material. We prove that two  \emph{boundary measures}  
give a  \emph{linear } Volterra integral equation of the first kind for the unknown    kernel. Hence, with two measures, the identification of the kernel, which in principle is a nonlinear problem,  is reduced to the solution of a mildly ill posed but \emph{linear} problem.
 \medskip

\section{Introduction}

Let $ \zthe=\zthe(x,t) $, $t>0$ and $x\in\ZOMq$ (a bounded region with smooth boundary).
The equation with persistent memory ($\Delta$ is the laplacian)
\begin{equation}\ZLA{eq:system}
\zthe'=\intt N(t-s)\Delta \zthe(s)\ZD s
 \end{equation}
 is encountered in thermodynamics of systems with memory   ($ \zthe  $ is the temperature)   in    nonfickian diffusion ($\zthe$ is the concentration)  or in viscoelasticity. In this case  $ \zthe(t) $ represents the displacement and the equation is usually written as
   \begin{equation}
\ZLA{eq:systemDIFFERENZIATO}
\zthe''=N(0)\Delta\zthe(t)+\intt M(t-s)\Delta\zthe(s)\ZD s\,,\qquad M(t)=N'(t)\,.
\end{equation}

  The vector
 \[ 
 q(x,t)=\intt N(t-s)\nabla \zthe(x,s)\ZD s
  \]  
  is the (density of the) flux at time $ t $ and position $ x $
      in thermodynamics or nonfickian diffusion. In viscoelasticity,
$ q'(x,t) $ is the traction. Let $ n(x) $ be the exterior normal to $ \partial\ZOMq  $ at the point $x$. The products $ q(x,t)\cdot n(x) $ or $ q'(x,t)\cdot n(x) $, $ x\in \partial\ZOMq $,
represent respectively the flux throughout the boundary or the traction on the boundary. In both the cases, they are observable quantities.

Especially in viscoelasticity,  $N(t)$ is called the \emph{ relaxation  kernel.} In this paper we assume that it is a smooth function (see below for details) and \emph{in the applications we mentioned  above, it must be $N(0)>0$.}

It is a fact that the kernel $ N(t)  $  is a material property. Only its  overall properties can be deduced from general principles (see for example~\cite{giorgiTermoProperties}) and the specific function $N(t)$ has to be identified from experiments taken on the system. The ``experiments'' usually considered    consist  in the observation of the traction/flux on the boundary of a suitable sample of the material, subject to known excitations. The memory kernel $ N(t) $ is identified from such observations. 

\begin{Remark}
{\rm The integrals in~(\ref{eq:system}) should extend from the time the material has been ``created'', conceivably $ -\infty $. We assume that the measures  are taken after a certain time $ t_0 $   sufficiently far  from the ``creation'' of the material and that the material was kept undisturbed for a long time, so that the effect of the memory of the previous history is negligible. This time $ t_0 $ we rename $ 0 $.

}
\end{Remark}

The shape of sample of  the material usually does  not influence the memory kernel      otherwise, $ N $ would depend also  on the position $ x $ (see~\cite{BykovMILLISECONDS}) so that we can try to determine $ N(t) $ from samples which have simple geometries. As in~\cite{BykovMILLISECONDS} and most of the references (commented in Sect.~\ref{sect:references}), we shall use samples which are in the shape of a slab, so that we can assume that $ \ZOMq $ is a segment and we can normalize
\[ 
\ZOMq=(0,\pi)\,.
 \] 
The identification of $ N(t) $ is a nonlinear problem since if $ N(t) $ is not yet known, the 
function $ \zthe (t) $ cannot be computed and the problem consists in the identification of the pair $ (N(t),\zthe(x,t)) $. 
 This pair enters in a nonlinear way in the equation. In spite of this, we shall see that using two measures we can construct a \emph{ linear} Volterra integral equation solely for $ N(t) $. So, \emph{the identification of $ N(t) $ is reduced to a linear problem.} This is a first difference with existing literature. A second difference is as follows: in many papers, especially  in engineering applications, the identification of $N(t)$ is done as follows: it is assumed that the kernel belongs to a certain class of functions which depend  on  a finite number of parameters.   
   The special kernel $N(t)$ assumed to be the ``true'' kernel is the one which minimizes the discrepancy between   a finite number of experimental data and the corresponding theoretically computed values.
 In contrast with this,   we don't assume any known a priori information on the functional form of $N(t)$, a part smoothness and the condition $N(0)>0$.

Due to different applications and terminology, for consistency and clarity we state that $\zthe$
will be called ``temperature'' and $q$ will be called the ``flux'' (note that the usual minus sign in the definition of the flux or the traction has no influence in the following arguments and it is ignored).

 \subsection{\ZLA{sect:assumpInfoDESCRI}The assumption on $ N(t) $ and informal description of the procedure}  
 
The first of our assumptions is that  the kernel $ N(t) $ is of class $ C^3 $ (so that the memory kernel $ M(t)=N'(t) $ of the viscoelastic system is of class $ C^2 $  as often assumed in this kind of problems, see~\cite{FabrizioOWENS}).
 
It is seen from the form~(\ref{eq:systemDIFFERENZIATO}) of Eq.~(\ref{eq:system}) that  this equation  is a perturbation of the wave equation \emph{since we explicitly assume} $ c^2=N(0)>0 $. A consequence of $ N(0)>0 $ is that \emph{signals propagate in the body with finite velocity $ c=\sqrt{N(0)} $.} So, whether this assumption holds   is easily checked and the value of $ c $ is easily computed: just measure the velocity of propagation of the signals in the sample of the material.

Our identification procedure applies in the case when $ c>0 $. Solely to present simpler formulas, we assume that a first, simple,  step  is the identification of the value   $ c =\sqrt{N(0)}$. We assume that this step has already been done. Then,  the change of the  time scale $ t\mapsto  t/c $ transform $ N(t) $ to  the new kernel $ N_1(t)=(1/c)N(t/c) $. Note that $ N_1(0)=1 $. So, we assume:
\begin{itemize}
\item the velocity $c>0$ of signal propagation has been identified.  
\item the transformation of the time scale has been done already and we rename as $ N(t) $ the function  $ N_1(t)= (1/c)N(t/c) $. So, we study Eq.~(\ref{eq:system}) (i.e.~(\ref{eq:systemDIFFERENZIATO})) with $ N(0)=1 $.

\end{itemize}
 
 Now we   describe informally the identification procedure.   We recall the goal: we fix any time $ T $ and we want to identify $ N(t) $   for $ t\in [0,T] $ (we already know $N(0)=1$).
 
\begin{enumerate}

\item We perform the following two measures, both on intervals of duration $ T $ (it is not restrictive to rename such intervals as $ [0,T] $):
 \begin{enumerate}
 
 \item  the temperature $ \zthe(x,t) $   at the ends $ x=0 $ and $ x=\pi $ is kept equal zero,  $\zthe(0,t)=0$, $\zthe(\pi,t)=0$.   We impose a \emph{special} initial condition $ \zthe(x,0)=\xi (x) $ and we read the flux $ y^{\xi } (t)= q(\pi,t) $ (note the lower case $y$). We use a suitable initial condition $\xi_0$ (a ramp) in order to identify  an intermediate function $K(t)$ which is then used in the identification of the memory kernel.
 
 \item  We impose null initial condition and  $ \zthe(\pi,t)=0 $ but now we excite the system from the left end: we impose $ \zthe (0,t)=f(t) $. We measure $ q(\pi,t) $. For clarity, the flux due to the boundary displacement $ f  $ (with zero initial condition) is denoted $ Y(t)=Y^f(t) $ (note the upper case).
 
 In this second measure we don't need to use any special function $f(t)$, in particular $f(t)$ is not related to the initial condition used in the first measure. But in practice the function $f(t)$ will be smooth. For continuity we have to take $f(0)=0$ and $f(t)$ ``saturates'' in a short time, i.e. it  approaches $\lim _{t\to+\ZIN} f(t)$ (which in practice is finite) quite fast.

\end{enumerate}
\item We derive an algorithm for the identification of $ N(t) $ using these two observations.
\end{enumerate}

  We note that  the first  step might be done by imposing $ \zthe(x,t)=0 $ and $ \zthe'(x,0)=\eta(x) $, with a special initial velocity $ \eta(x) $ (integrating~(\ref{eq:systemDIFFERENZIATO}) we see that in this case Eq.~(\ref{eq:system})  must have an additive term $ \eta(x) $ on the right hand side). For definiteness, we study $ \eta=0 $ and   $ \xi \neq 0$  but in practical applications to viscoelasticity $ \xi=0 $ and $ \eta\neq 0 $ might be easier to impose.

Now we describe the organization of the paper. In Sect.~\ref{sect:references} we discuss few references on the identification problem and in Sect.~\ref{Sect:PRELIMINARIES} we present the preliminary
 properties of system~(\ref{eq:system}) we shall needed. The identification algorithm is in Sect.~\ref{Sect:IdeProCE}.
 
 It is possible to identify $N(t)$ using a similar algorithm, which uses two measures with different initial conditions (and  zero temperature at the boundary). We believe that this algorithm is not practical since  to impose initial condition   is more difficult  then to impose a boundary temperature (or deformation). However, for completeness, this second algorithm is described in Sect.~\ref{sect:differentALGO}.

\subsection{\ZLA{sect:references}Comments on previous references}

Due to the importance of this problem, reference on kernel identification is enormous and both Engineers and Mathematicians addressed the problem, along different and scarcely related lines. Most of the available results can be described as follows. An approach to kernel identification, seemingly first proposed 
in~\cite{LorenziSINESTRARI} and then followed by a wealth of related papers  (see also the books~\cite{FabrizioOWENS,Lerenzilibro} for references)   is as follows: it is noted that the identification of $N(t)$ on the basis of additional measures is in fact the identification of the pair $(\zthe(t),N(t))$ and $\zthe(t)$ solves an integrodifferential equation. Now the trick is to compute derivatives of the measured output (of course when the data are sufficiently smooth) and to derive a second integrodifferential equation for $N(t)$. In this way we get a   \emph{nonlinear} integrodifferential equation for the pair   $(\zthe(t),N(t))$ and, in general using fixed point theorems, it is possible to prove   existence of solutions of this systems.   

This leads to the identification in particular of $N(t)$ (according to the assumptions, either locally, for a short time, or globally, as in~\cite{LoreSine2}).

The identification requires only one measure, taken assuming known (and smooth) initial condition (or affine term).

The approach that is most common in engineering literature is as follows. It is assumed that the kernel $N(t)$ belongs to a certain class of functions, identified by a finite number of parameters. The classes that are most often considered are Dirichlet sums, i.e.
\begin{equation}\ZLA{eq:Diri}
N(t)=\sum_{k=0}^{K} a_k e^{-b_k t}
\end{equation}
as for example in~\cite{BykovMILLISECONDS,Gerlac,Sellier,ZhuDapengHONEYcomb}
(here $a\geq 0$, $b\geq 0$) or a combination of Abel kernels
 \begin{equation}
 \ZLA{eq:abelKern}
 N(t)=\sum_{k=0}^{K}
 \frac{a_k }{t^{\beta_k}} 
 \end{equation}
sometimes a combination of both  (see for example~\cite{Dinzart,GolubKozbarRegulina}).   

The class of the kernels~(\ref{eq:Diri}) or~(\ref{eq:abelKern})  have a ``physical'' sense since they are encountered when solving ordinary differential equations of integer or fractional order and  these classes of kernels are reported to fit sufficiently well experimental data, but not in every application, see~\cite{Ruymbeke}.

In order to identify the parameters, a sample of the material, usually in the form of a slab, is excited 
  using   a known initial or boundary deformation. The    stress relaxation  is measured on (a part of the) boundary and the  parameters ($a_k$, $b_k$, sometime also the exponents $\beta_k$) are chosen so to fit better the experimental measures (usually minimizing a quadratic index). This can be done using a single measure. Multiple measures can be used to improve accuracy or to remove spurious frequencies.  
  
  Different algorithms works in the frequency domain, using Fourier or Laplace transformations, as for example in~\cite{Janno}.
  
  The use of a single measure is common to the algorithm presented in~\cite{LorenziSINESTRARI,LoreSine2} and related papers. Here instead we assume that \emph{two} different measure are performed and, from these measures, the kernel $N(t)$ is reconstructed. We don't assume that the kernel belongs to any specific class of functions, a part regularity and the physical property $N(0)>0$. In fact, no further properties of the kernel are in principle needed for the reconstruction algorithms (not even positivity for every time!) in spite  of the fact that such conditions are imposed by thermodynamics and can be used to improve the results, see
  Sect.~\ref{sect:simulation}. \emph{Note however that we are assuming that $N(t)$ is smooth for $t\geq 0$ and this rules out Abel kernels.} The extension of our algorithm to Abel kernels will be pursued in future researches.

\section{\ZLA{Sect:PRELIMINARIES}Preliminary properties of system~(\ref{eq:system})}

We present a representation of the solutions of system~(\ref{eq:system}) when $\ZOMq=(0,\pi)$, taken from~\cite{PandIEOT,Pandcina,LoretiPANDOLFIsforza} (see~\cite{PandolfiLIBRO,PandSHARP} for more general results, in particular for extensions to the case ${\rm dim}\, \ZOMq>1$).
\bigskip

We consider the sequence of the functions

\[
\phi_n=\sqrt{\frac{2}{\pi}}\sin nx\,,\qquad n\in \mathbb{N}\,.
\]
It is known form the theory of Fourier series that $\{\phi_n\}$ is an orthonormal basis of $L^2(0,\pi)$ and that
\[
\Delta \phi_n(x)=\frac{\ZD^2}{\ZD x^2} \phi_n(x)=-n^2 \phi_n(x)\,.
\]
We expand the initial condition $\xi$ and the solution in sine series:
\[
\xi_n=\intp \xi(x)\sin nx\ZD x\,,\qquad \zthe_n(t)=\intp \zthe(x,t)\sin nx\ZD x
\]
\[
\frac{\pi}{2}\xi(x)=\ZSUno \xi_n\sin nx\,,\qquad \frac{\pi}{2}\zthe(x,t)=\ZSUno \zthe_n(t)\sin nx
\]
\[
\zthe'_n=-n^2\intt N(t-s) \zthe_n(s)\ZD s+n\intt N(t-s) f(s)\ZD s\,,\qquad \zthe_n(0)=\xi_n\,.
\]
This is a scalar integrodifferential equation whose solution can be represented as

\[
 \zthe_n(t)=  \xi_n z_n(t) +  \left [ n \intt z_n(t-s)\ints N(s-r) f(r)\ZD r\,\ZD s\right ]
\]
where $z_n(t)$ solves
\begin{equation}\ZLA{equaDIz}
z_n'(t)=-n^2\intt N(t-s) z_n(s)\ZD s\,,\qquad z_n(0)=1\,.
\end{equation}

The transformation $(\xi,f)\mapsto w$ is linear and continuous from  $L^2(0,\pi)\times L^2(0,T)$ to $C([0,T]; L^2(0,\pi))\times C^1([0,T];H^{-1}(0,\pi))$.

We need the following observation:
\begin{Lemma} Let $T>0$. There exists $M=M_T$ such that for every $t\in[0,T]$ and every $n$ we have:
 
 \begin{equation}\ZLA{DiseFONDAdaPROVARE}
\left |\frac{z_n'(t)}{n} +e^{-\zaa t}\sin nt\right |\leq \frac{M}{n}\,,\qquad \zaa=-\frac{1}{2}N'(0)\,. 
 \end{equation}
\end{Lemma}
\zProof
 This inequality follows for example from the equality~(5.4) in~\cite{PandIEOT} (see also~\cite{Pandcina,LoretiPANDOLFIsforza}). We must be careful with the notation. With $\zaa=-N'(0)/2$ let $\tilde N(t)=e^{2\zaa t}N(t)$ so that $\tilde N(0)=1$ and $\tilde N'(0)=0$. The functions denoted $z_n(t)$ and $N(t)$ in~\cite[formula~(5.4)]{PandIEOT} are respectively $e^{2\zaa t} z_n(t)$ and $\tilde N(t)=e^{2\zaa t}N(t)$ (see~\cite[formula~(5.1)]{PandIEOT}). So,  we have the following equality for the function $z_n(t)$ in~(\ref{equaDIz}) (here    $\zben=\sqrt{n^2-\zaa^2}$):
   \begin{align*}
 z_n(t)&=  e^{-\zaa t}\cos\zben t+\frac{\zaa}{\zben}e^{-\zaa t} \sin\zben t 
  -\frac{n^2}{\zben^2} \intt e^{-2\zaa s } z_n(s)\left [
 \tilde N'(t-s)+\right.\\
 &\left. +\int_0^{t-s}\left [\zaa \tilde N'(t-s-r)-\tilde N''(t-s-r)\right ]e^{\zaa r} \cos\zben r\ZD r
 \right ]\ZD s \,.
   \end{align*}
   We compute the derivatives of both the sides, using $\tilde N'(0)=0$. We get the equality
   \begin{align*}
   z_n'(t)+ \zben e^{-\zaa t}\sin \zben t &= -\frac{\zaa^2}{\zben} e^{-\zaa t} \sin\zben t -\\&
   -\frac{n^2}{\zben} \intt e^{-2\zaa s} z_n(s)\left \{
  \tilde N''(t-s) 
  -\tilde N''(0)e^{\zaa(t-s)}\cos\zben(t-s)\right.\\
 &\left.+\int_0^{t-s}\left[
  \zaa \tilde N''(t-s-r)-\tilde N'''(t-s-r)
  \right ] e^{\zaa r}\cos\zben r\ZD r
   \right\}\ZD s\,.
   \end{align*}
   
   The result follows by dividing both the sides  with $n$ and nothing that if $t\in [0,T]$ then there exists $M$ such that
   \[
|\sin\zben t-\sin t|\leq \frac{M}{n}\,.   \zdiaform
   \]

We are going to use these properties of the system for the identification of the kernel $N(t)$, assumed unknown. Note that these same  properties have also been used to solve a problem of source identification, see~\cite{PANDsurvey,PandDCDS2}.

\section{\ZLA{Sect:IdeProCE}The identification procedure}

Now we give the details of the two steps of the identification procedure described in Sect.~\ref{sect:assumpInfoDESCRI}.

\subsection{\ZLA{Sect:ideTWOiniCOND}
A first measure of the flux, using  initial temperature}
 
Using a special initial temperature, we can identify an intermediate function $K(t)$ which will be used in the second step of kernel identification.

Let $\zthe(x,0)= \xi(x)=\frac{2}{\pi}\ZSUno \xi_n\sin n x \in L^2(0,\pi) $ so that $ \{\xi_n\}\in l^2 $ and
 
\begin{align*} 
\frac{\pi}{2}\zthe(x,t)&= \ZSUno \xi_n z_n(t)\sin nx
\in C([0,T];L^2(0,\pi))
\,,\\ 
\frac{\pi}{2}\zthe_x(x,t)&=\ZSUno (n\xi_n)z_n(t)\cos nx\in C([0,T];H^{-1}(0,\pi))\,.
 \end{align*}
 So we have
 \begin{equation}\ZLA{eq:flussoDAxiANTE} 
\frac{\pi}{2} q(x,t)=\ZSUno  \xi_n \left [n\intt N(t-s)z_n(s)\ZD s\right ]\cos nx=-\ZSUno\frac{1}{n} \xi_n z'_n(t)\cos nx  
  \end{equation}
  and
  \begin{equation}\ZLA{eq:flussoDAxi}
  y^\xi(t)=\frac{\pi}{2}q(\pi,t)=-\ZSUno \frac{1}{n}(-1)^n\xi_n z_n'(t)\,.  
  \end{equation}
Using the inequality~(\ref{DiseFONDAdaPROVARE}), we see   that $y^\xi(t)\in L^2(0,T)$ for every  $T>0$ and every $\xi\in L^2(0,\pi)$.  
  If in particular we choose
  \begin{equation}\ZLA{ramp}
\xi_0(x)=\frac{1}{2}(\pi-x)=\ZSUno\frac{1}{n}\sin nx  
  \end{equation}
  then   the series in~(\ref{eq:flussoDAxi}) is
  \[
  -\ZSUno \frac{1}{n^2}(-1)^nz_n'(t)\,.
  \]
  The function we wanted to identify is
  \begin{equation}\ZLA{definizK}
K(t)=\ZSUno (-1)^n\frac{1}{n^2}z_n'(t)=-y^{\xi_0}(t) \,.  
  \end{equation}
  
  This is the first measure we need in the identification process.

 \subsection{\ZLA{Sect:ALGOritmFINALE}Measure of the flux due to a boundary excitation}
 
Now we combine the identification of $K(t)$ obtained  in the previous section   with an observation of the flux at $ x=\pi $ due to a nonzero temperature applied to the left hand $ x=0 $ (and zero initial condition). Hence we consider the following system:
\begin{equation}
\ZLA{eq:perFLUssoDEFORMA}
\left\{\begin{array}{l}
\displaystyle   \zthe'(t)=\intt N(t-s)\Delta\zthe(s)\ZD s\,,\quad \zthe(x,0)=0\,,\\
\displaystyle  \zthe(0,t)=f(t)\,,\quad \zthe(\pi,t)=0\,.
\end{array}\right.
\end{equation}
We choose any $ C^1 $ input $ f $ such that $ f(0)=0 $, (this  is consistent with the fact that the initial condition is zero). So we have
\[ 
f(t)=\intt g(s)\ZD s
 \]
 and we explicitly assume
 \[ g(0)=f'(0)\neq 0\,. \]
 We use 
 \begin{align}
\nonumber \frac{\pi}{2} \zthe(x,t)&=\ZSUno (\sin nx)\intt f(s)\left [n\int_0^{t-s}N(r)z_n(t-s-r)\ZD r\right ]\ZD s=\\
 \nonumber &=\ZSUno (\sin nx)\left [\frac{1}{n}\intt f(s)\frac{\ZD}{\ZD s} z_n(t-s)\ZD s\right ]=\\
\ZLA{eq:DefiSOLuzW} &= \left( \ZSUno\frac{1}{n} (\sin nx)\right )\left (\intt g(r)\ZD r-\intt g(s)z_n(t-s)\ZD s\right )\,.
  \end{align}
  The function $ (x,t)\mapsto  \zthe(x,t) $ belongs to $ C([0,T];L^2(0,\pi) )$
 so that $ (x,t)\mapsto  \zthe_x(x,t) $ belongs to $ C([0,T];H^{-1}(\ZOMq) )$,
 \begin{align*}
\frac{\pi}{2}  \zthe_x(x,t)&=\left (\intt g(r)\ZD r\right )\ZSUno \cos nx-\ZSUno \left (\cos nx\right )\intt g(s) z_n(t-s)\ZD s=\\
 &=\pi\ZDE(x)\left (\intt g(r)\ZD r\right )-\frac{1}{2}\left (\intt g(r)\ZD r\right ) -
 \ZSUno \left (\cos nx\right )\intt g(s) z_n(t-s)\ZD s\,.
  \end{align*}
 Here $ \ZDE(x) $ is the Dirac's delta (supported at $ x=0 $) and so when observing the flux at $x=\pi$ we get
 \begin{align}
\nonumber& \frac{\pi}{2}Y^f(t)= \frac{\pi}{2}q(\pi,t) =-\frac{1}{2}\left (\intt N(t-s)\ints  g(r)\ZD r\,\ZD s\right )-\\
\nonumber &-\ZSUno(-1)^n \intt g(s)\int_0^{t-s}N(t-s-r)z_n(r)\ZD r\,\ZD s= \\
\nonumber  &=-\frac{1}{2}\left (\intt N(t-s)\ints g(r)\ZD r\,\ZD s\right )+\ZSUno (-1)^n\frac{1}{n^2} \intt g(s)z_n'(t-s)\ZD s=\\
\nonumber&=-\frac{1}{2}\left (\intt N(t-s)\ints g(r)\ZD r\,\ZD s\right )+  \intt g(t-s)\left [\ZSUno (-1)^n\frac{1}{n^2}z_n'( s)\right ]\ZD s\\
\ZLA{equa:significatoK}  &= \intt g(t-s)\left ( K(s)-\frac{1}{2}\ints N(r)\ZD r\right )\ZD s 
  \end{align}
  where, from~(\ref{definizK}),
    \[ 
  K(s)=\left [\ZSUno (-1)^n\frac{1}{n^2}z_n'( s)\right ]\,.
   \]
   We have seen in the previous section   that $ K(t) $ can be estimated from the measure of one boundary flux. 
   So, at this point, \emph{$ K(t) $ is a known function.}   
Formula~(\ref{equa:significatoK}) explains its significance: $K(t)-1/2\intt N(s)\ZD s$ is the impulse response of the system with boundary input $ f $ at $ x=0 $ and the flux observed at $ x=\pi $.

 \begin{Remark}
 {\rm
 In fact, the $\ZDE$-function is an artifact, due to the fact that formula~(\ref{eq:DefiSOLuzW}) gives an odd extension of the solution which is periodic of period $ 2\pi $ and so the formula introduces a jump at $ t=0 $ (because $ f(t)\neq 0 $) and this produces a $ \ZDE $ function at $ x=0 $ and also at every multiple of $ 2\pi $. This fact has been discussed in~\cite{AvdonPAND1,AvdonPand2}.\zdia
 }
 \end{Remark}

The function $ N(t) $ can be identified from formula~(\ref{equa:significatoK})
which is a Volterra integral equation of the first kind for 
\[
M(t)=\intt N(s)\ZD s\,.
\]
Although this may not be the best way from the point of view of numerical computation, we may note that if $ g\in C^1$, hence $ f\in C^2 $, then $ Y(t) $ is   differentiable and  upon computing  the derivatives of both the sides of~(\ref{equa:significatoK}) we get:
\begin{align*}
g(0) M(t)+\intt g'(t-s) M(s)\ZD s=\\
 =
 2 g(0)K(t)+2 \intt g'(t-s)K(s)\ZD s
-\pi\left (\frac{\ZD }{\ZD t} Y^f(t) \right 
 )\,.
 \end{align*}
 We have chosen an input signal $f$ such that $f'(0)=g(0)\neq 0$. So,
this is a Volterra integral equation of the second kind (but note that the right hand side requires the numerical computation of $Y'(t)$).

Once   this Volterra integral equation has been solved, a last step of numerical differentiation identifies $ N(t) $.
\begin{Remark}
{\rm We note that
\begin{itemize}
\item
If $g(0)=1$ and $g'(t)=0$ then we get an explicit formula for $M(t)$ but in this case $f(t)=t$. A steadily increasing temperature cannot be applied. In practice, $f(t)$ has to saturate, as for example $f(t)=t/(t+1)$ or $1-e^{-t}$.
\item In practice it does not seems difficult to impose a variable temperature (or traction) on the boundary. It may be more difficult to impose a special initial condition but we note that the initial condition we need in the identification of $K$ is a ramp, which may be achieved by letting the equilibrium be reached between two different temperatures of  the ends of the sample.\zdia
\end{itemize}
}
\end{Remark}

   \subsection{\ZLA{sect:differentALGO}A different identification procedure}
   
   We present now a second identification algorithm which uses two different measures, obtained  using two different, but related, initial conditions.
  In the algorithm we proposed above, the most questionable step seems  the use of a special initial condition.   It is even more  difficult to impose also a second, related, initial temperature.     So, the procedure we present now may not have  a practical interest, but it seems mathematically interesting.

   We impose an initial temperature $\xi(x)=\ZSUno \xi_n\sin nx$ and we use~(\ref{eq:flussoDAxi}):
 \begin{align*}
 y^\xi(t)&=  
-\ZSUno (-1)^n
\xi_n\frac{z_n'(t)}{n}
  \,.
  \end{align*}
  Let us compute the convolution $ N* y^\xi$ (we recall that the kernel $ N(t) $ is still unidentified):
\begin{equation} 
\ZLA{eq:CONVOLdatINIZ}
 \intt N(t-s)y^\xi(s)\ZD s =
 -\ZSUno \frac{(-1)^n}{n} \xi_n\intt N(s)z_n'(t-s)\ZD s 
\,.
\end{equation}
This formula is easily interpreted: let us consider the second initial condition 
\[ 
\eta(x)=  \ZSUno \frac{ 	\xi_n}{n^2}\sin nx =x\left (\ZSUno\frac{\xi_n}{n}\right )-\int_0^x\int_0^s \xi(r)\ZD r\,\ZD s
 \]
 and let us compute $y^{\eta}(t)$ from the intermediate equality in~(\ref{eq:flussoDAxiANTE}) (where we put $x=\pi$ and we replace $\xi_n$ with  $\xi_n/n^2$). We get
 \[ 
 y^\eta(t)=\ZSUno (-1)^n \frac{\xi_n}{n} \intt N(s) z_n(t-s)\ZD s
  \]
  so that 
  \begin{align}
  \frac{\ZD}{\ZD t}y^\eta(t)&=\left (\ZSUno\frac{(-1)^n \xi_n}{n}\right )N(t)+\ZSUno(-1)^n\frac{\xi_n}{n} \intt N(s)z_n'(t-s)\ZD s\\
&  = \left ( \ZSUno\frac{(-1)^n\xi_n}{n}\right )N(t)-\intt y^\xi(t-s)N(s)\ZD s\,.
   \end{align}

   The functions $ y^\xi $ and $ y^\eta $ can be  measured and so  we find a Volterra integral equation for the unknown function $ N(t) $. The Volterra integral equation is of the second kind, i.e. its solution is a well posed problem, if $ \xi $ is such that
   \[ 
   \ZSUno\frac{(-1)^n \xi_n}{n}\neq 0\,.
    \]
   In spite of this, the identification problem is ill posed since the procedure requires the numerical computation of the derivative of $y^\eta(t)$.

\subsection{\ZLA{sect:simulation}A toy simulation}
Whether the identification procedure we presented is really useful is a matter to be settled by the Engineers on the basis of concrete experiments. Here we present a toy  simulation  on the basis of theoretical computations.
  We assume that $ N(t) $ is a linear combination of exponentials, as it is often assumed in applications, and for simplicity we confine ourselves to the case that $ N(t) $ is the combination of three exponentials,
 
 \begin{equation}\ZLA{Cha2:KerneESEMP}
 N(t)=c_1 e^{- \zt_1 t}+c_2 e^{-\zt_2 t}+c_3e^{-\zt_3 t} 
  \end{equation}
  where the coefficients $c_k$ and the exponents $\zt_k$ are positive.
  Then, the equation of $ z_n(t) $ is
  \[ 
  z_n'(t)=-n^2\intt \left [ c_1 e^{-\zt_1(t-s)}+c_2e^{-\zt_2(t-s)}+c_3 e^{-\zt_3(t-s)}  \right ] z_n(s)\ZD s\,.
   \]
It is easy to reduce the computation of $ z_n(t) $ to the solution of an ordinary differential equation. In fact,   
 $z_n(t)$ is the first component of the solution $(z,w_1,w_2,w_3)  $ of the following problem:
  \[ 
    \left\{\begin{array}
{l}
 z '=-n^2(c_1 w_1+c_2w_2+c_3w_3) \\
w_1'=-\zt_1w_1+z   \\
w_2'=-\zt_2 w_2+z    \\
w_3'=-\zt_3w_3+z   
\end{array}\right. 
\qquad  
\left(\begin{array}{c}
z(0)\\w_1(0)\\w_2(0)\\w_3(0)
\end{array}\right)=
\left(\begin{array}{c}
1\\0\\0\\0
\end{array}\right)\;.
   \]
   So, the functions $z_n(t)$ can be easily computed and   we can compute also $Y^f(t)$  and $K(t)$.   
   Then we choose $f(t)= 1-e^{-t}$ and we reconstruct the kernels $M(t)=\intt N(s)\ZD s$  and $N(t)$.
 
 Before we present the simulation we note that:
 \begin{itemize}
 \item the reconstruction of the relaxation kernel is in two steps: first we identify $M(t) $ and then a next step of numerical differentiation gives $N(t)$. The identification of $M(t)$ is a deconvolution problem which can be solved, essentially with the same results, with Tikonov regularization (as described in~\cite{Groetch}) or Lavrientev method (see~\cite{FagnPAND1,FagnPAND2}). Both these methods introduce fast oscillations  but it is known in this   problem that the oscillations    are an artifact due to the amplification of the noise. In fact, $M(t)$ is increasing (since for physical reasons $N(t)\geq 0$) and, in the physical case that $N(t)$ is decreasing, $M(t)$ is concave. Hence these oscillations can be attenuated for example by computing the average of the reconstruction of $M(t)$ on nearby steps. This has been done in the simulation below.
 \item in the simulation we present below we choose to reconstruct $M(t)$ using Lavrentev method (followed by an average, to remove fast oscillations). Instead we use Tikonov regularization in the next step of the reconstruction of $N(t)$ from  $M(t)$   since Lavrentev method is more efficient when it is known that the function to be reconstructed is $0$ at the initial time (and since Lavrentev method is more sensitive to disturbances).
 \item Also the function $N(t)$ can be averaged to reduce the effect of errors since in practice we know that it is a decreasing (and often a convex) function. Taking into account the known value $N(0)$ this last step (not shown in the plots) will also improve the reconstruction near $t=0$.
 \end{itemize}
 
\begin{figure}[h]
\caption{\label{figCASO1}Left: the kernel $M(t)$ (the function, the reconstruction and its average); right: the reconstruction of the relaxation kernel $N(t) $ (and its plot from the expression~(\ref{eq:IlNUcleoExpo}) } \vspace{0cm}
\begin{center}
\includegraphics[width=6cm]{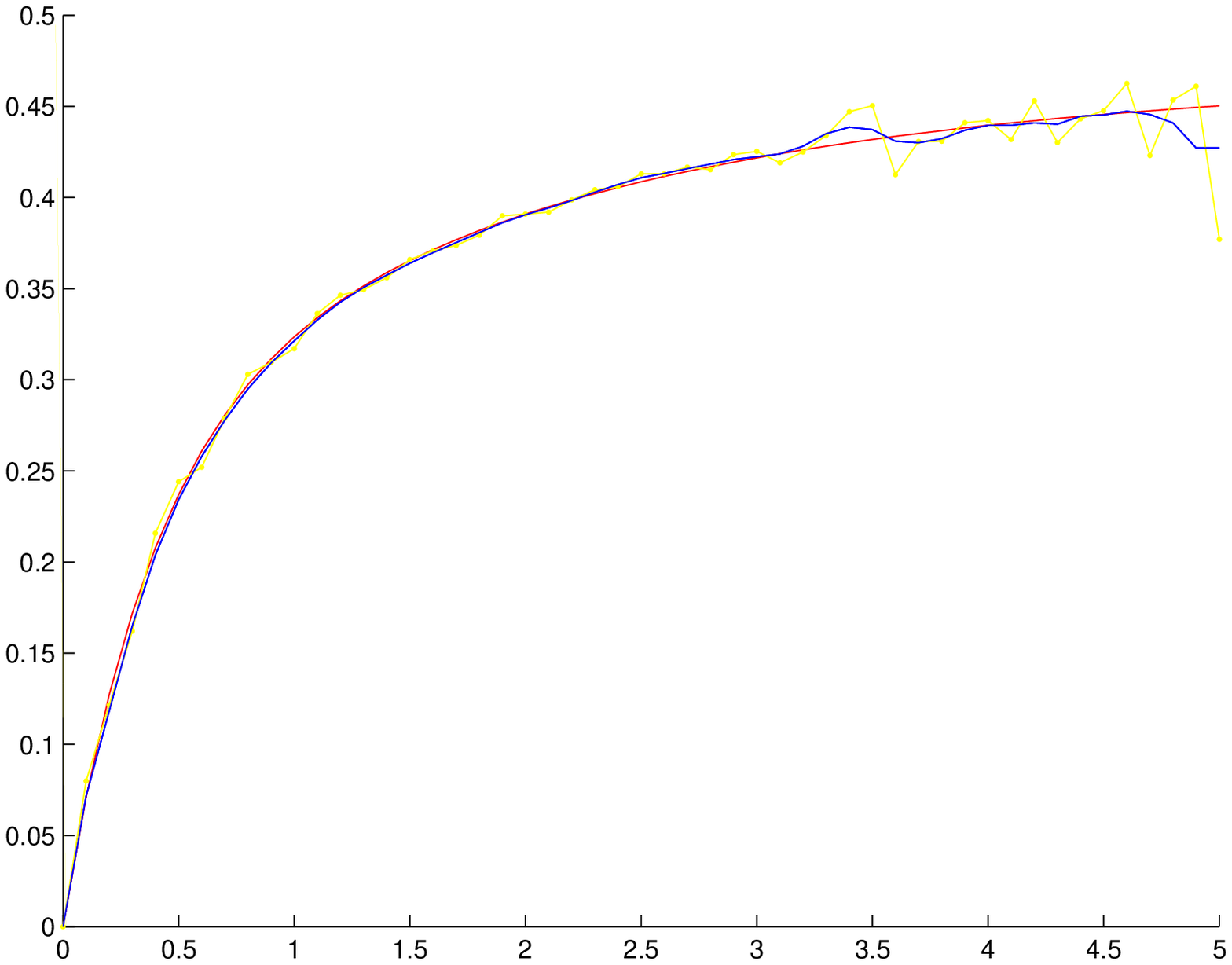} 
\includegraphics[width=6cm]{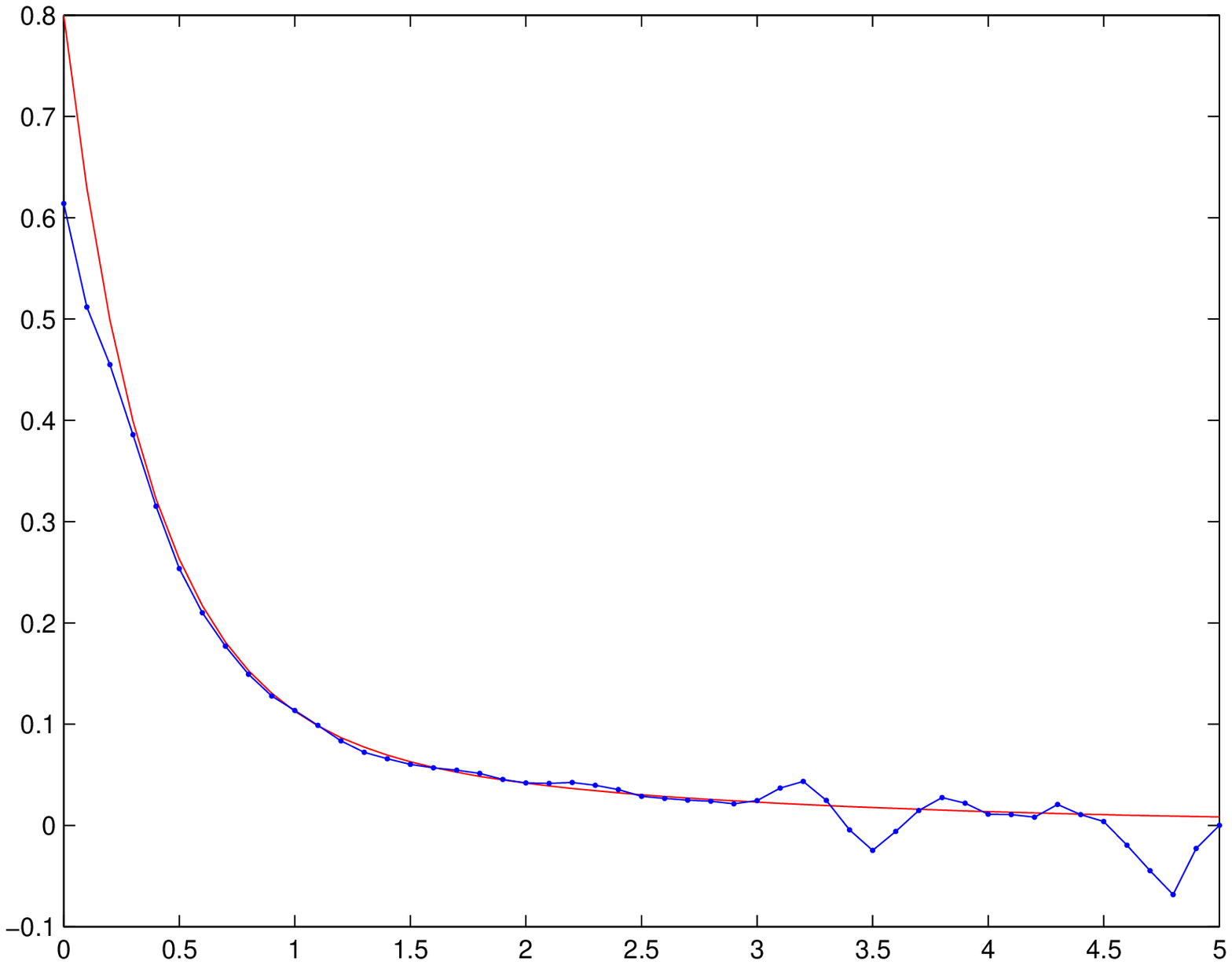}
\end{center}
%\vspace{-1cm}
\smallskip
\end{figure}

 The special case we test is as follows: the time interval is $[0,5]$ and we assume $10$ measures every time unit.  
 We choose 
 \begin{equation} \ZLA{eq:IlNUcleoExpo}
N(t)= (1/10)e^{-t/2}+(1/5)e^{-2t}+(1/2)e^{-3t} 
 \end{equation}
 (i.e. we normalize  so to have $N(0)=1$). The functions $K(t)$ and $Y^f(t)$ are approximated (and corrupted by  an error of $1\%$ at every ``measure''). Then, we ``reconstruct'' both $M(t)$ and $N(t)$. 
   Figure~\ref{figCASO1} (left) presents the plots of $M(t)$ and its reconstruction from Lavrientev method (the plot where the amplification of the disturbances is evident)  and its average (at every step tree values on the left and on the right have been averaged). The averaged reconstruction of $M(t)$ is fed to the Tikonov algorithm in order to compute the numerical derivative of $M(t)$, i.e. $N(t)$, which is plotted on the right. The penalization parameters of the Lavrientev and Tikonov algoritms are respectively $0.1$ (equal to the interval of time between consecutive measures) and $0.01$.
   
 \enddocument